\documentclass[a4paper,12pt]{article}
\usepackage{amsthm}
\usepackage{amssymb}
\usepackage{amsmath}
\usepackage{amsfonts}

\newtheorem{theorem}{Theorem}[section]

\newtheorem{corollary}{Corollary}[section]

\def\b1{\mbox{\boldmath $1$}}

\newenvironment{demo*}{\vspace{3mm}\noindent{\bf Proof.}}{\hfill $\Box$ \vspace{3mm}}

\begin{document}

\title{\bf Asymptotic results for tail probabilities of sums
of dependent heavy-tailed random variables}

\author{\normalsize$^a$Kam Chuen Yuen,  and $^b$Chuancun Yin\\
%\footnote{Corresponding author: E-mail address: ccyin@mail.qfnu.edu.cn}\\
{\normalsize\it $^a$Department of Statistics and Actuarial Science,
The University of Hong Kong,} \\
\noindent{\normalsize\it Pokfulam Road, Hong Kong}\\
{\normalsize\it E-mail: kcyuen@hku.hk} \\
[3mm]
{\normalsize\it $^b$School of Mathematical Sciences, Qufu Normal University}\\
{\normalsize\it Shandong 273165, P.R.\ China}\\
{\normalsize\it E-mail: ccyin@mail.qfnu.edu.cn }}
\date{}
\maketitle
%\centerline {\bf Running title:  Asymptotics for tail probabilities}

\vskip0.6cm
 \noindent{\large {\bf Abstract}} \ Let $\{X_1, X_2,
\cdots \}$ be a sequence of dependent heavy-tailed random variables
with distributions $F_1, F_2,\cdots$ on $(-\infty,\infty)$, and let
$\tau$ be a nonnegative integer-valued random variable independent
of the sequence $\{X_k, k \ge 1\}$. In this framework, we study the
asymptotic behavior of the tail probabilities of the quantities
$X_{(n)} = \max_{1\le k \le n} X_k$, $S_n =\sum_{k=1}^n X_k$ and
$S_{(n)}=\max_{1\le k\le n} S_k$ for $n>1$, and for those of their
randomized versions $X_{(\tau)}$, $S_{\tau}$ and $S_{(\tau)}$.   We
also consider applications of the results obtained to  some
commonly-used risk processes.



%\vfill \noindent \footnotesize {${}^*$ Corresponding author: e-mail
%address: ccyin\@mail.qfnu.edu.cn \\
%\phantom{*} Mathematics Subject Classification (2000) 62E20 $\cdot $
%62H20 $\cdot $ 62P05}

\newpage
\normalsize

\baselineskip=20pt

%\noindent{\bf 1.~~Introduction}
\section{Introduction}\label{intro}

Throughout this paper, let $\{X_1, X_2, \cdots\}$ be a sequence of
random variables with distributions $F_1, F_2,\cdots$ supported  on
$\mathbb{R}:=(-\infty,\infty)$ satisfying
$\overline{F}_k(x)=1-F_k(x) > 0$ for all $x$.  For $n\ge 1$, we
write
\[
X_{(n)}=\max_{1\le k\le n}X_k, \quad S_n = \sum_{k=1}^{n}X_k, \quad
S_{(n)}=\max_{1\le k\le n}S_k.
\]
Let $\tau$ be a counting random variable independent of $\{X_k, k\ge
1\}$. Then, the randomized versions of $X_{(n)}$, $S_n$, $S_{(n)}$
are given by $X_{(\tau)}$, $S_{\tau}$, $S_{(\tau)}$. Tail
probabilities of the quantities $X_{(n)}$, $S_n$, $S_{(n)}$,
$X_{(\tau)}$, $S_{\tau}$, $S_{(\tau)}$ with heavy-tailed random
variables are of great interest in finance, insurance and many other
disciplines. Since accurate distributions for these quantities are
not available in most cases, the study of asymptotic relationships
for their tail probabilities becomes important. Many results have
been derived under different degrees of generality in the
literature. In particular, most of the results are for independent
$X_1, \cdots, X_n$ with distributions belonging to the class of
subexponential distributions.

For two independent random variables $X$ and $Y$ with distribution
functions $F$ and $G$ supported on $(-\infty,\infty)$, we denote by
$F*G(x) =\int_{-\infty}^{\infty}F(x-y)dG(y)$, $-\infty<x<\infty$,
the convolution of $F$ and $G$, and by $F^{*n}=F*F\cdots *F$ the
$n$-fold convolution of $F$. Unless otherwise stated, all limit
relations are  for $x\to\infty$.  By definition, a distribution $F$
on $[0,\infty)$ is said to be subexponential ($F\in \cal {S}$) if
the relation $\overline{F^{\ast 2}}(x)\sim 2 \overline{F}(x)$
($x\to\infty$) holds where the symbol $\sim$ means that the ratio of
the two sides tends to 1. More generally, a distribution function
$F$ on $(-\infty,\infty)$ belongs to the subexponential class
$\cal{S}$ if $F^+(x)=F(x)1(x\ge 0)$ does, where $1(\cdot)$ is the
indicator function. A recent account on tail asymptotic results for
the sum of two independent random variables can be found in Foss and
Korshunov [11]. They proved that, for two distributions $F_1$ and
$F_2$ on $[0,\infty)$, if one of $F_1$ and $F_2$ is heavy tailed
(see the definition below), then
\begin{equation}
\liminf_{x\to\infty}\frac{\overline{F_1*F_2}(x)}{\overline{F_1}(x)
+\overline{F_2}(x)}=1, \label{intro-eq1}
\end{equation}
and that, for any heavy-tailed distribution $F$,
\begin{equation}
\liminf_{x\to\infty}\frac{\overline{F^{*2}}(x)}{\overline{F}(x)}=2.
\label{intro-eq2}
\end{equation}
Denisov et al. [6] extended (\ref{intro-eq2}) to
\begin{equation}
\liminf_{x\to\infty}\frac{\overline{F^{*\tau}}(x)}{\overline{F}(x)}=E
\tau, \label{intro-eq3}
\end{equation}
with $\tau$ being a light-tailed random variable. Furthermore, for
any heavy-tailed distribution $F$ on $\mathbb{R}^+$ with finite
mean, Denisov et al. [7] showed that if $P(c\tau>x)=o(\overline
{F}(x))$ for some $c>EX$ as $x\to\infty$, then (\ref{intro-eq3})
holds. Also, if $F$ is subexponential and $\tau$ is light tailed and
independent of the summands, then
\begin{equation}
P(S_{\tau}>x)\sim \overline{F}(x)E\tau, \;
x\to\infty.\label{intro-eq4}
\end{equation}
Note that all the above-mentioned results were established for
independent nonnegative random variables. When the random variables
are possibly negative and dependent according to certain structures,
the validity of these results remains to be studied.

We now consider three examples in which some of the above relations
do not hold. The first comes from Yu et al. [25] while the last two
are extracted from Albrecher et al. [1].

\noindent{\bf Example 1.1.} \ Assume that $X$ is a discrete random
variable with masses $p_n=(X=2^{n+1}-1)=2^{-n-1}$, $n\ge 0$. Denote
its distribution by $\rho$. Then, for any $0<q<1$, define
 $F=q\rho+(1-q)\sigma$, where $\sigma$ is a non-degenerated distribution on a subset of $(-\infty,0]$.
  Without loss of generality, we assume that $\sigma$ has support on $[-3,0)$, and that $\sigma (-2)-\sigma(-3)=\delta>0$. Then, $F$ is heavy-tailed but does not belong to the class $\mathcal{L}$ of distributions with long tails ($F\notin \mathcal{L}$). Also, it satisfies
 $$\liminf_{x\to\infty}\frac{\overline{F^{*2}}(x)}{\overline{F}(x)}<2.$$

\noindent{\bf Example 1.2.} \ Let $X_1$ and $X_2$ have marginal
distribution function $F$ belonging to the subexponential class
$\mathcal{S}$.
 Then, there exists a copula for $X_1$ and $X_2$  such that
 $$ \lim_{x\to\infty}\frac{P(X_1+X_2>x)}{P(X_1>x)}=\infty.$$

\noindent{\bf Example 1.3.} \ Assume that random variables $X_1$ and
$X_2$  are comonotone dependent with common distribution $F =
1-\overline{F}$ where $\overline{F}\in\mathcal{R}_{-\alpha}$ is
regularly varying at infinity with some index $\alpha>0$. Then, we
have $$\lim_{x\to\infty}\frac{P(X_1+X_2>x)}{P(X_1>x)}=2^{\alpha}.$$

\noindent These examples indicate that relations (1.1)-(1.4) may not
hold for heavy-tailed distributions supported on $[a,\infty)$ with
$a < 0$ or for dependent heavy-tailed distributions.  More examples
can be found in Albrecher et al. [1].

The purpose of this paper is to find sufficient conditions under
which relations (\ref{intro-eq1})-(\ref{intro-eq4}) hold for
possibly negative, non-identically distributed, and dependent
heavy-tailed random variables. The outline of the paper is as
follows. Section 2 presents several classes of heavy-tailed
distributions and the dependence assumptions used in later sections.
Section 3 is devoted to the tail behavior of $X_{(n)}$, $S_n$ and
$S_{(n)}$. Section 4 investigates the tail behavior of $X_{(\tau)}$,
$S_{\tau}$ and $S_{(\tau)}$. Section 5 presents applications of the
main results to  risk theory.

%\vskip 0.1cm
\section{Preliminaries}\label{prelimin}
\setcounter{equation}{0}

A random variable $X$ (or its distribution $F$) is heavy tailed (to
the right) if $E\exp(\alpha X)=\infty$ for all $\alpha>0$, and light
tailed otherwise.
 One of the most important classes of heavy-tailed
distributions is the subexponential class. Closely related to the
class $\mathcal{S}$ are the class $\mathcal{S}^*$, the class
$\mathcal {D}$ of distributions with dominatedly varying tails, and
the class $\mathcal{L}$ of distributions with long tails. By
definition, a distribution function $F$ on $\mathbb{R}$ with finite
mean belongs to the class $\mathcal{S}^*$ if and only if
$\overline{F}(x) > 0$ for all $x$ and $\int_0^x
\overline{F}(x-y)\overline{F}(y)dy\sim 2m_{{F}^+}\overline{F}(x),$
as $x\to\infty,$ where $m_{{F}^+}=\int_0^{\infty}\overline{F}(x)dx$
is the mean of $F^+$. It is known that if $F\in \mathcal{S}^*$, then
both $F$ and $F_I$ are subexponential, where $F_I$ is defined by
$\overline{F_I}(x)=\min(1,\int_x^{\infty}\overline{F}(t)dt)$; see
Kl\"uppelberg [15]. A distribution function $F$ with support on
$(-\infty,\infty)$ belongs to the class $\mathcal{D}$ if
$\limsup_{x\to\infty}{\overline{F}(xy)}/{\overline{F}(x)}<\infty$
holds for some (or equivalently for all) $0<y<1$. Obviously, if
$F\in\mathcal{D}$, then, for any $y>0$, $\overline{F}(xy)$ and
$\overline{F}(x)$ are of the same order as $x\to\infty$ in the sense
that
$$0<\liminf_{x\to\infty}\frac{\overline{F}(xy)}{\overline{F}(x)}\le
\limsup_{x\to\infty}\frac{\overline{F}(xy)}{\overline{F}(x)}<\infty.$$
A distribution function $F$ is said to belong to the class
$\mathcal{L}$ if
$\lim_{x\to\infty}{\overline{F}(x+y)}/{\overline{F}(x)}=1$ holds for
some (or equivalently for all) $y$. One can easily check that for a
distribution $F\in \mathcal{L}$, there exists a positive function
$h(x)\to\infty$ such that $\overline{F}(x+h(x))\sim
\overline{F}(x).$ The class $\mathcal{S}^*$ and the intersection
$\mathcal{D}\cap\mathcal{L}$ are two well-known subclasses of
subexponential distribution functions. For details of these classes
of heavy-tailed distributions and their applications, the reader is
referred to Asmussen [2],  Bingham et al. [3], Embrechts et al. [9],
and Embrechts et al. [10].
 Furthermore, a distribution $F$ is said to
be strongly subexponential, denoted by $F\in \mathcal{S}_*$, if
$\overline{F_h^{*2}}(x) \sim 2 \overline{F_h}(x),$ uniformly in
$h\in [1,\infty)$, where the distribution $F_h$ is defined as
\begin{equation}
F_h(x)=\min\left(1, \int^{x+h}_x F(t)dt\right), \quad  x>0.
\label{ee}
\end{equation}
See Korshunov [17] for sufficient conditions for some distribution
to belong to the class $\mathcal{S}_*$. Kaas and Tang [14] proved
that $\mathcal{S}_*$ is a subclass of $\mathcal{S}$ while Denisov et
al. [5] showed that $\mathcal{S}^*$ is a subclass of
$\mathcal{S}_*$.  The relations between the above-mentioned classes
are summarized as follows:
$\mathcal{D}\cap\mathcal{L}\subset\mathcal{S}_*\subset\mathcal{S}\subset\mathcal{K}
$, and that, if the distribution function $F$ has a finite mean,
then $F\in\mathcal{D}\cap\mathcal{L}\Rightarrow F\in \mathcal{S}^*
\subset\mathcal{S}_*,$ where $\mathcal{K}$ represents the class of
distribution functions with heavy tails.

Recall that $X_1,\cdots, X_n$ are $n$ real-valued random variables
with distributions $F_1,\cdots, F_n$, respectively.   Here, we
assume that these random variables are dependent.
 To model the dependence of a multivariate
distribution with non-identical marginals, one may use the theory of
copulas (see, for example, Nelsen [19]). A copula is a multivariate
joint distribution defined on the $n$-dimensional unit cube
[0,1]$^n$ such that every marginal distribution is uniform on the
interval [0,1]. By Sklar's theorem, for a  multivariate joint
distribution $F$ of a random vector $(X_1,\cdots,X_n)$ with
marginals $F_1, \cdots, F_n$, there exists a copula $C$ such that
\begin{equation}
F(x_1, \cdots, x_n)=C(F_1(x_1),\cdots,F_n (x_n)), \quad (x_1,\cdots,
x_n)\in \mathbb{R}^n. \label{copula}
\end{equation}
If $F_1, \cdots, F_n$  are all continuous, then $C$ is unique and
can be written as
\[
C(u_1,\cdots,u_n)=P(F_1(X_1)\le u_1, \cdots, F_n(X_n)\le
u_n)=F(F_1^{-1}(u_1),\cdots, F_n^{-1}(u_n)),
\]
for any $(u_1,\cdots, u_n)\in [0,1]^n$. Conversely, if $C$ is a
copula and $F_1, \cdots, F_n$ are distribution functions, then $F$
defined in (\ref{copula}) is a multivariate joint distribution with
marginals $F_1, \cdots, F_n$.

For notional convenience, we state the following four  assumptions
regarding the  random variables $X_1,\cdots, X_n$.

\noindent {\bf H1.} \ Assume that $X_1,\cdots, X_n$ satisfy
\[
\lim_{x\to\infty}\frac {P(X_i>x,
X_j>x)}{\overline{F_i}(x)+\overline{F_j}(x)}=0,
\]
for all $1\le i\neq j\le n$.    \hfill $\Box$

This dependence assumption was   first introduced in Chen and Yuen
[4] and    $\{X_i\}$  is called pairwise quasi-asymptotically
independence.

 \noindent {\bf
H2.} \ Assume that
\begin{equation}
\hat{\lambda}_{ij}=\lim_{x_i\wedge x_j\to\infty} P(|X_i|>x_i
|X_j>x_j)=0. \label{prelimin-eq1}
\end{equation}
holds for all $1\le i\neq j\le n$, where $x_{i}\wedge
x_{j}=\min(x_{i}, x_{j})$. This concept is related to the so-called
asymptotic independence; for example, see Resnick [23].  Note that
the asymptotic independence means a large value in one component is
unlikely to be accompanied by a large value in another. \hfill
$\Box$

\noindent {\bf H3.} \ Assume that there exist positive constants
$x_0$ and $c_0$ such that the inequality
\[
P(X_i>x_i|X_j=x_j, j\in J)\le c_0\overline {F_i}(x_i)
\]
holds for all $1\le i\le n, \emptyset \neq J\subset \{1,2,\cdots,
n\}\backslash \{i\}, x_i>x_0$, and $x_j >x_0$ with $j\in J$. When
$x_j$ is not a possible value of $X_j$, the conditional probability
above  is simply understood as 0.  Note that this dependence
assumption was used in Geluk and Tang [14]. \hfill $\Box$

\noindent {\bf H4.} \ Let $X_1,\cdots, X_n$ be dependent. Assume
that the dependent structure is governed by an absolutely continuous
copula $C(u_1,\cdots,u_n)$ such that there exists  positive constant
$M<\infty$ with $c(u_1,\cdots,u_n)\le M$ for all
$(u_1,\cdots,u_n)\in [0,1]^n$, where $c$ is the copula density given
by
\[
c(u_1,\cdots,u_n)=\frac{\partial ^n C(u_1,\cdots, u_n)}{\partial
u_1\cdots\partial u_n}.
\]
\hfill $\Box$

\medskip
\noindent {\bf Remark 2.1.} \ It is obvious that ${\bf H2}$ implies
${\bf H1}$. Also, we see in Geluk and Tang [13] that ${\bf H3}$
implies ${\bf H2}$.  If $J$ is the set of a single point, then one
can show that ${\bf H4}$ implies ${\bf H3}$. In fact, let
$(X_1^*,\cdots, X_n^*)$ be an independent copy of $(X_1,\cdots,
X_n)$, then
\begin{eqnarray*}
P(X_i>x_i|X_j=x_j, j\in J)&=&\lim_{\Delta x_j\to 0}P(X_i>x_i|x_j\le
X_j<x_j+\Delta x_j, j\in J)\\
&=&\lim_{\Delta x_j\to 0}\frac{P(X_i>x_i, x_j\le X_j<x_j+\Delta x_j,
j\in J)}{P(x_j\le X_j<x_j+\Delta x_j, j\in J)}\\
&\le & M \lim_{\Delta x_j\to 0}\frac{P(X^*_i>x_i, x_j\le
X^*_j<x_j+\Delta x_j,
j\in J)}{P(x_j\le X_j<x_j+\Delta x_j, j\in J)}\\
&\le & M P(X^*_i>x_i),
\end{eqnarray*}
as desired. Note that from the proof above, it is  easy to see that
if additionally $c$ satisfies $c(u_1,\cdots,u_n)\ge m$ for all
$(u_1,\cdots,u_n)\in [0,1]^n$, where $m$ is a positive constant,
then  ${\bf H4}$ implies ${\bf H3}$.\hfill $\Box$

To end the section, we present an example in which the four
assumptions are satisfied.

\noindent {\bf Example 2.1.} \ A joint $n$-dimensional distribution
is called a Farlie-Gumbel-Morgenstern (FGM) distribution if it has
the form
\begin{equation}
F(x_1,\cdots,x_n)=C(F_1(x_1),\cdots,F_n(x_n))), \quad (x_1,\cdots,
x_n)\in \mathbb{R}^n, \label{prelimin-eq2}
\end{equation}
where $F_1,\cdots, F_n $ are the one-dimensional marginals, and the
copula $C$ is given by
\begin{equation}
C(u_1,\cdots,u_n)=\prod_{k=1}^n u_k \left(1+\sum_{1\le i<j\le
n}a_{ij}(1-u_i)(1-u_j)\right), \quad (u_1,\cdots, u_n)\in [0,1]^n,
\label{prelimin-eq3}
\end{equation}
where $a_{ij}$ are real numbers fulfilling certain requirements so
that $F(x_1,\cdots,x_n)$ is a proper $n$-dimensional distribution.
For details of FGM distributions, see Kotz et al. [18].

It is easy to check that if the random variables $X_1,\cdots, X_n$
follow a joint $n$-dimensional FGM distribution defined in
(\ref{prelimin-eq2}) and (\ref{prelimin-eq3}) whose marginal
distributions $F_k$ ($1\le k\le n$) are absolutely continuous and
satisfy $F_k(-x)=o(\overline{F_k}(x))$, then the four assumptions
{\bf H1-H4} are fulfilled. \hfill $\Box$

\section{Results for finite sums}\label{finitely}
\setcounter{equation}{0}

In this section, we present our main results for finite sums.

Let $(X_1^*,\cdots, X_n^*)$ be an independent copy of $(X_1,\cdots,
X_n)$, that is, $(X_1^*,\cdots, X_n^*)$ and $(X_1,\cdots, X_n)$ are
two independent random vectors with the same marginal distributions
and the components of $(X_1^*,\cdots, X_n^*)$ are independent.
Similar to $S_n$ and $S_{n,k}$, we define $S_{n}^*= \sum_{k=1}^n
X_k^*$, and $S_{n,k}^*=S_n^*-X_k^*$.

\begin{theorem}\label{thrm3-1} \ Let $X_1, \cdots, X_n$ be heavy-tailed random variables with
distribution $F_1,\cdots,$ $ F_n$, respectively. Under  assumption {\bf H3}, \\
(1) if  $F_1, F_2,\cdots, F_n\in \cal{L}$, then
\begin{equation}
\liminf_{x\to\infty}\frac{P(S_n>x)}{\sum_{k=1}^n \overline{F_k}(x)}=
\liminf_{x\to\infty}\frac{P(S_{(n)}>x)}{\sum_{k=1}^n
\overline{F_k}(x)}=1; \label{finitely-eq1}
\end{equation}
(2) if additionally  $P(\sum_{k=1}^n
X_k^{*+}>x)\sim\sum_{k=1}^{n}\overline{F_k}(x)$, then
\begin{equation}
\lim_{x\to\infty}\frac{P(S_n>x)}{\sum_{k=1}^n \overline{F_k}(x)}=
\lim_{x\to\infty}\frac{P(S_{(n)}>x)}{\sum_{k=1}^n
\overline{F_k}(x)}=1. \label{finitely-eq2}
\end{equation}
\end{theorem}

\noindent {\bf Proof.} \ We first show that if  $F_1, F_2,\cdots,
F_n\in \cal{L}$, then
\begin{equation}
\liminf_{x\to\infty}\frac{P(S_n>x)}{\sum_{k=1}^n
\overline{F_k}(x)}\ge 1, \label{lemma3-1-eq1}
\end{equation}
under assumption {\bf H2}.

It follows from the definition of the class  $\cal{L}$ that there
exists a sequence $\{a(x)\}$ such that  $a(x)\to\infty$ as
$x\to\infty$, $2a(x)\le x$, and
\[
\overline{F}_k(x+a(x))\sim \overline{F_k}(x),\quad x\to\infty, \quad
k=1,2,\cdots, n.
\]
Note that
\begin{eqnarray*}
P(S_n>x)
& \ge    & P(S_n>x, X_{(n)}>x+a(x))\\
& \ge    & \sum_{k=1}^n P(S_n>x, X_k>x+a(x))\\
&        & \quad - \ \sum_{1\le i<j\le n}P(X_i>x+a(x), X_j>x+a(x))\\
& \equiv & I_1 (x)+I_2 (x).
\end{eqnarray*}
Assumption {\bf H2} implies that
\[
I_2 (x)=o\left(\sum_{k=1}^n\overline{F_k}(x)\right).
\]
Write $S_{n,k}=S_n-X_k$ for $1\le k\le n$. Then,
\begin{eqnarray*}
I_1 (x)
& \ge & \sum_{k=1}^n P(S_{n,k}>-a(x), X_k>x+a(x))\\
& =   & \sum_{k=1}^n P(X_k>x+a(x))- \sum_{k=1}^n P(S_{n,k}\le -a(x),
X_k>x+a(x)).
\end{eqnarray*}
It follows from (\ref{prelimin-eq1}) that
\[
\sum_{k=1}^n P(S_{n,k}\le -a(x),
X_k>x+a(x))=o\left(\sum_{k=1}^n\overline{F_k}(x+a(x))\right).
\]
Thus,
\begin{eqnarray*}
I_1 (x)
& \ge  & \sum_{k=1}^n P(X_k>x+a(x))-o\left(\sum_{k=1}^n\overline{F_k}(x+a(x))\right)\\
& \sim & \sum_{k=1}^n\overline{F_k}(x).
\end{eqnarray*}
This proves (\ref{lemma3-1-eq1}).  Note that (\ref{lemma3-1-eq1})
also holds under assumption {\bf H3} as {\bf H3} implies {\bf H2}.

We next show that if $F_1, F_2,\cdots, F_n\in \cal{L}$ and
assumption {\bf H3} holds, then
\begin{equation}
\liminf_{x\to\infty}\frac{P(S_n>x)}{\sum_{k=1}^n
\overline{F_k}(x)}\le 1. \label{lemma3-3-eq1}
\end{equation}
If additionally  $P(\sum_{k=1}^n
X_k^{*+}>x)\sim\sum_{k=1}^{n}\overline{F_k}(x)$, then
\begin{equation}
\limsup_{x\to\infty}\frac{P(S_n>x)}{\sum_{k=1}^n
\overline{F_k}(x)}\le 1. \label{lemma3-3-eq2}
\end{equation}

Assume that at least one of  $X_1,\cdots, X_n$ has infinite mean. In
this special case,
 (\ref{lemma3-3-eq1}) still holds without assumption {\bf H3}. In fact,
 when $X_1,\cdots, X_n$ are nonnegative, one  obtains
\[
\frac{\int_0^t P(S_n>x)dx}{\sum_{k=1}^{\infty}\int_0^t
\overline{F_k}(x)dx} = \frac{E\min(S_n,t)}{\sum_{k=1}^{\infty}E
\min(X_k,t)}\le 1,
\]
for any positive $t$. This together with $\int_0^{\infty}
\overline{F_1}(y)dy=\infty$ imply that
\[
\limsup_{x\to\infty}\frac{P(S_n>x)}{\sum_{k=1}^n
\overline{F_k}(x)}\le 1,
\]
and hence
\[
\liminf_{x\to\infty}\frac{P(S_n>x)}{\sum_{k=1}^n
\overline{F_k}(x)}\le 1.
\]

If at least one of  $X_1^+ = X_1 1(X_1\ge 0),\cdots, X_n^+ = X_n
1(X_n\ge 0)$ has infinite mean, we may consider  $X_1^+, \cdots,
X_n^+$. If at least one of $X_1,\cdots, X_n$ has negative infinite
mean, we may consider $-X_1,\cdots, -X_n$. So, in either case, we
have upper bound (\ref{lemma3-3-eq1}).

Now, suppose that all random variables  $X_1,\cdots, X_n$ have
finite means. It is clear that the inequality
\begin{equation}
P(X^+_{(n)}>x)\le \sum_{k=1}^n \overline{F_k}(x).
\label{lemma3-3-eq3}
\end{equation}
holds. The conditions $F_1, F_2,\cdots, F_n\in \cal{L}$ and {\bf H3}
imply that
 there exist positive constants $x_0$ and $d_n$ such that
$$P(S_{n,k}>x|X_k=x_x)\le d_n P(S^*_{n,k}>x)$$
holds for all $1\le k\le n, x>x_0$ and $x_k>x_0$; see Geluk and Tang
[13]. It can be shown that   for every function $a(\cdot) :
(0,\infty)\rightarrow (x_0,\infty)$ and for every $1\le k\le n$, and
$x>x_0$,
\begin{equation}
P\left(S_n >x,a(x)<X_k \le x\right)\le d_n P\left(S_n^*
>x, a(x)<X^{*}_k \le x\right). \label{lemma3-3-eq4}
\end{equation}
Using (\ref{intro-eq1})  and induction arguments, it follows from
the proof of Lemma~5.2 of Geluk and Tang [13] that there exists a
sequence $\{x_{l}\}$ such that $l\to\infty$ and
\begin{equation}
P\left(S_n^*>x_{l}, \frac{x_{l}}{n}<X_j^*\le
x_{l}\right)=o\left(\sum_{k=1}^n \overline{F_k}(x_{l})\right).
\label{lemma3-3-eq5}
\end{equation}
In particular, if \ $P(S_n^*>x)\sim \sum_{k=1}^n \overline{F_k}(x)$,
then
\[
P\left(S_n^*>x, \frac{x}{n}<X_j^*\le x\right)=o\left(\sum_{k=1}^n
\overline{F_k}(x)\right).
\]
Hence, (\ref{lemma3-3-eq3}), (\ref{lemma3-3-eq4}) and
(\ref{lemma3-3-eq5}) imply  that
\begin{eqnarray*}
P(S_n>x_{l})
& \le  & P\left(\cup_{k=1}^n (X_k^+>x_{l})\right)+P\left(S_n^+>x_{l}, \cap_{k=1}^n (X_k^+\le x_{l})\right)\\
& \le  & P\left(X^+_{(n)}>x_{l}\right)+ \sum_{k=1}^n P\left(S_n^+>x_{l}, \frac{x_{l}}{n}<X_k^+\le x_{l}\right)\\
& \le  & P\left(X^+_{(n)}>x_{l}\right)+ d_n\sum_{k=1}^n P\left(S^{+*}_n>x_{l}, \frac{x_{l}}{n}<X^{+*}_k\le x_{l}\right)\\
& \sim & \sum_{k=1}^n \overline{F_k}(x_{l}).
\end{eqnarray*}
This proves (\ref{lemma3-3-eq1}), and (\ref{lemma3-3-eq2}) can be
proved using similar arguments.

Finally, (\ref{finitely-eq1}) and (\ref{finitely-eq2}) follow from
(\ref{lemma3-1-eq1}), (\ref{lemma3-3-eq1}), (\ref{lemma3-3-eq2}) and
the fact that
\[
P\left(\sum_{k=1}^n X_k^+>x\right)\ge P(S_{(n)}>x)\ge P(S_n>x).
\]
This completes the proof of Theorem~\ref{thrm3-1}.  \hfill $\Box$

\medskip
\noindent {\bf Remark 3.1.} \ Let $X_1,\cdots, X_n$ be $n$
nonnegative heavy-tailed random variables.  It follows from Foss and
Korshunov [11] that there exist two sequences $\{x_{l}\}$ and
$\{a(x_{l})\}$ such that $x_{l}\to\infty$, $a(x_{l})\to\infty$ as
$l\to\infty$, $2a(x_{l})\le x_{l}$, and
\[
\overline{F}_k(x_{l}+a(x_{l}))\sim \overline{F_k}(x_{l}),\quad
l\to\infty, \quad k=1,2,\cdots, n.
\] Repeating the same  arguments as above we find that (3.3) also
holds without the assumption that $F_1, F_2,\cdots, F_n\in \cal{L}$.

Note that for $n=2$ and nonnegative random variables, a result
similar to (\ref{finitely-eq2}) was established by Foss and
Korshunov [11] under the independence assumption.  \hfill $\Box$

\begin{corollary}\label{corollary3-1} \
Let $X_1,\cdots, X_n$ be heavy-tailed random variables with common
distribution $F$. Under assumption {\bf H3}, if  $F\in \cal{L}$,
then
\begin{equation}
\liminf_{x\to\infty}\frac{P(S_n>x)}{\overline{F}(x)}=\liminf_{x\to\infty}\frac{P(S_{(n)}>x)}{\overline{F}(x)}=n;
\label{finitely-eq3}
\end{equation}
and if $F\in\cal{S}$, then
\begin{equation}
\lim_{x\to\infty}\frac{P(S_n>x)}{\overline{F}(x)}=\lim_{x\to\infty}\frac{P(S_{(n)}>x)}{\overline{F}(x)}=n.
\label{finitely-eq4}
\end{equation}
\end{corollary}

If we switch our attention from the class $\mathcal{L}$ to the class
$\mathcal{D}$, we only need assumption {\bf H2} which is weaker than
{\bf H3}.

\begin{theorem}\label{thrm3-2} \
Let $X_1,\cdots, X_n$ be $n$ nonnegative heavy-tailed random
variables with distribution $F_1,\cdots,F_n$, respectively. Under
assumption {\bf H2}, if $F_k \in{\mathcal{D}}$ for $k = 1,\cdots,n$,
then
\begin{equation}
\liminf_{x\to\infty}\frac{P(S_n>x)}{\sum_{k=1}^n
\overline{F_k}(x)}=\liminf_{x\to\infty}\frac{P(S_{(n)}>x)}{\sum_{k=1}^n
\overline{F_k}(x)}=1; \label{th3-eq1}
\end{equation}
Moreover, if $F_k \in{\mathcal{L}}$ for $k = 1,\cdots,n$,
 then
\begin{equation}
\lim_{x\to\infty}\frac{P(S_n>x)}{\sum_{k=1}^n
\overline{F_k}(x)}=\lim_{x\to\infty}\frac{P(S_{(n)}>x)}{\sum_{k=1}^n
\overline{F_k}(x)}=1. \label{th3-eq2}
\end{equation}
\end{theorem}
\noindent {\bf Proof.} \ It follows from Foss and Korshunov [11]
that there exist two sequences $\{x_{l}\}$ and $\{a(x_{l})\}$ such
that $x_{l}\to\infty$, $a(x_{l})\to\infty$ as $l\to\infty$,
$2a(x_{l})\le x_{l}$,  and
\begin{equation}\overline{F}_k(x_{l}-a(x_{l}))\sim
\overline{F_k}(x_{l}),\;k=1,2,\cdots, n. \label{th3-eq3}
\end{equation}
Note that
\begin{eqnarray*}
P(S_n>x_{l})
& \le    & P(S_n>x_{l}, X_{(n)}>x_{l}-a(x_{l}))+P(S_n>x_{l}, \cap_{k=1}^n (X_k\le x_{l}-a(x_{l}))\\
& \le    & \sum_{k=1}^n P(X_k>x_{l}-a(x_{l})) +\sum_{k=1}^n P(S_n-X_k>a(x_{l}), X_k>\frac{x_{l}}{n})\\
& \equiv & J_1 (x_{l})+J_2 (x_{l}).
\end{eqnarray*}
It is obvious that $J_1 (x_{l})\sim \sum_{k=1}^n
\overline{F_k}(x_{l}).$ Recall that $S_{n,k}=S_n-X_k$ for $1\le k\le
n$. Under assumption {\bf H2}, we have
\begin{eqnarray*}
J_2 (x_{l})
& =   & \sum_{k=1}^n P\left(S_{n,k}>a(x_{l}), X_k>\frac{x_{l}}{n}\right)\\
& \le & \sum_{1\le k\neq j\le n}P\left(X_j>\frac{a(x_{l})}{n-1}, X_k>\frac{x_{l}}{n}\right)\\
& \le & \sum_{1\le k\neq j\le
n}P\left(X_j>\frac{a(x_{l})}{n-1}|X_k>\frac{x_{l}}{n}\right)
P\left(X_k>\frac{x_{l}}{n}\right)\\
& =   & o(1)\sum_{k=1}^n \overline{F_k}(x_{l}),
\end{eqnarray*}
where,  in the last step, we used the relation
$\overline{F}_k(xy)=O(\overline{F}_k(x))$ for all $y>0$ (as  $F_k
\in{\mathcal{D}}$). Thus, we obtain
\[
\liminf_{x\to\infty}\frac{P(S_n>x)}{\sum_{k=1}^n
\overline{F_k}(x)}\le 1.
\]
This together with (\ref{lemma3-1-eq1}) (see Remark 3.1) imply that
\[
\liminf_{x\to\infty}\frac{P(S_n>x)}{\sum_{k=1}^n
\overline{F_k}(x)}=1.
\]
Furthermore, if $F_k \in{\mathcal{L}}$ for $k = 1,\cdots,n$, then
(\ref{th3-eq3}) holds for any $x$. Repeating the steps above, we
obtain
\[
\limsup_{x\to\infty}\frac{P(S_n>x)}{\sum_{k=1}^n
\overline{F_k}(x)}\le 1.
\]
The remaining proofs of (\ref{th3-eq1}) and (\ref{th3-eq2}) are
straightforward. \hfill $\Box$

The last result in this  section is trivial.

\begin{theorem}\label{thrm3-3} \
Let $X_1,\cdots, X_n$ be heavy-tailed random variables with
distribution $F_1,\cdots,$ $ F_n$, respectively. Under assumption
{\bf H1}, we have
$$\lim_{x\to\infty}\frac{P(X_{(n)}>x)}{\sum_{k=1}^n
\overline{F_k}(x)}=1.$$
\end{theorem}
%\medskip
\noindent {\bf Proof.} \ On one hand, the inequality
$P(X_{(n)}>x)\le \sum_{k=1}^n \overline{F_k}(x)$ is trivial. On the
other hand,
\[
P(X_{(n)}>x)\ge\sum_{k=1}^n P(X_k>x)-\sum_{1\le i\neq j\le
n}P(X_i>x,X_j>x).
\]
Under assumption  {\bf H1}
\[
\sum_{1\le i\neq j\le n}P(X_i>x,X_j>x) \le o(1)\sum_{k=1}^n
\overline{F_k}(x).
\]
The result of Theorem~\ref{thrm3-3} follows from  the above
inequalities. \hfill $\Box$

\medskip
\noindent {\bf Remark 3.2.} \ For heavy-tailed random variables
$X_1, \cdots, X_n$ with distribution $F_1,\cdots,$ $ F_n$, if $F_1,
F_2,\cdots, F_n\in \cal{S}$, $F_i*F_j\in\cal{S}$ for all $1\le i\neq
j\le n$, and assumption {\bf H2} holds, then the result of Geluk and
Tang ([13, Theorem 3.2]) gives
$$\lim_{x\to\infty}\frac{P(S_n>x)}{\sum_{k=1}^n \overline{F_k}(x)}=1.$$
Obviously, the conditions in Theorem \ref{thrm3-1} are slightly more
general than the above conditions of Geluk and Tang ([13, Theorem
3.2]). Furthermore,  (\ref{finitely-eq4})  proves  the insensitivity
of relation $P(S_n>x)\sim P(S_{(n)}>x)\sim n\overline{F}(x)$ to the
dependence assumption in {\bf H2}, and Theorem \ref{thrm3-3} proves
the insensitivity of relation $P(X_{(n)}>x)\sim n\overline{F}(x)$ to
the dependence assumption in {\bf H1}.     \hfill $\Box$

\setcounter{equation}{0}
\section{Results for random sums}\label{randomly}
Parallel to the result of Denisov et al. [6] for i.i.d. random
variables with a support on $[0,\infty)$, we obtain the following
theorem for heavy-tailed random variables satisfying some dependence
structure.

\begin{theorem}\label{thrm4-1} \
Let $X_1,\cdots, X_n$ be heavy-tailed random variables with common
distribution $F$, and $\tau$ be a counting random variable
independent of the sequence $\{X_k\}$. Assume that
$E(z^{\tau})<\infty$ for some $z>1$. Under assumption {\bf H3}, if
$F\in\cal{L}$, then
\begin{equation}
\liminf_{x\to\infty}\frac{P(S_{\tau}>x)}{\overline{F}(x)}=E\tau;\label{randomly-eq1}
\end{equation}
and if $F\in\cal{S}$, then
\begin{equation}
\lim_{x\to\infty}\frac{P(S_{\tau}>x)}{\overline{F}(x)}=E\tau.\label{randomly-eq2}
\end{equation}

\end{theorem}

\medskip
\noindent {\bf Proof.} \ Note that the condition on
$E(z^{\tau})<\infty$ implies the existence of $E\tau$ and
$E(1+\varepsilon)^{\tau}<\infty$ for some sufficiently small
$\varepsilon>0$. For such $\varepsilon$,  it follows from
(\ref{finitely-eq3}) that there exist a constant $K\equiv
K(\varepsilon)<\infty$ and a sequence $\{x_l\}$ such that for all
$n\ge 2$
\[
P(S_n>x_l)\sim n\overline{F}(x_l),\;\;P(S_n>x_l) \le
K(1+\varepsilon)^n\overline{F}(x_l),\;\; l\ge 1.
\]
Applying the dominated convergence theorem yields
\[
\lim_{l\to\infty}\frac{P(S_{\tau}>x_l)}{\overline{F}(x_l)}
=\lim_{l\to\infty}\sum_{n=1}^{\infty}P(\tau=n)\frac{P(S_{n}>x_l)}{\overline{F}(x_l)}=E\tau,
\]
which implies (\ref{randomly-eq1}). Furthermore, if $F\in\cal{S}$,
then (\ref{randomly-eq2})  follows from (\ref{finitely-eq4}) easily.
\hfill $\Box$

Here, we also extend the results of  Denisov et al. [7] to the case
of dependent heavy-tailed random variables. Note that their
asymptotic results are for sums of random size $\tau$ of i.i.d.
nonnegative heavy-tailed random variables with $\tau$ belonging to
the class of all light-tailed distributions and also to some class
of heavy-tailed distributions.

\begin{theorem}\label{thrm4-2} \
Let $X_1,\cdots, X_n$ be heavy-tailed random variables with common
distribution $F$, and $\tau$ be a counting random variable
independent of the sequence $\{X_k\}$ with infinite mean. Under
assumption {\bf H2}, we have
\begin{equation}
\liminf_{x\to\infty}\frac{P(S_{\tau}>x)}{\overline{F}(x)}=\liminf_{x\to\infty}\frac{P(S_{(\tau)}>x)}{\overline{F}(x)}
=\lim_{x\to\infty}\frac{P(X_{(\tau)}>x)}{\overline{F}(x)}=E\tau.
\label{randomly-eq3}
\end{equation}
\end{theorem}
%\medskip
\noindent {\bf Proof.} \ Since $\tau$ is independent of $X_k$'s, we
can write
\[
P(S_{\tau}>x)=\sum_{n=0}^{\infty}P(\tau=n)P(S_n>x).
\]
Using (\ref{lemma3-1-eq1}) and Fatou's lemma, we obtain
\[
\liminf_{x\to\infty}\frac{P(S_{\tau}>x)}{\overline{F}(x)}\ge
\sum_{n=0}^{\infty}n P(\tau=n)=E\tau.
\]
We immediately get the first equality in (\ref{randomly-eq3}) since
$E\tau=\infty$. The rest of the proof is similar. \hfill $\Box$

The next two theorems consider the case with $E\tau<\infty$.

\begin{theorem}\label{thrm4-3} \
Let $X_1,\cdots, X_n$ be heavy-tailed random variables with common
distribution $F$, and $\tau$ be a counting random variable
independent of the sequence $\{X_k\}$. Assume that $E\tau<\infty$.
Under assumption {\bf H1}, we have
\begin{equation}
\lim_{x\to\infty}\frac{P(X_{(\tau)}>x)}{\overline{F}(x)}=E\tau.\label{randomly-eq4}
\end{equation}
\end{theorem}
\noindent {\bf Proof.} \ The proof is straightforward.
  \hfill $\Box$

\begin{theorem}\label{thrm4-4} \
Let $X_1,\cdots, X_n$ be heavy-tailed random variables with common
distribution $F$ and finite mean, and $\tau$ be a counting random
variable  independent of the sequence $\{X_k\}$ with finite mean
$E\tau$. Under assumption {\bf H4}, we have the following results:
\begin{enumerate}
\item[\rm {(i)}]  Assume that $EX_1<0$. If $F\in\cal{S}$, then
\begin{equation}
\liminf_{x\to\infty}\frac{P(S_{\tau}>x)}{\overline{F}(x)}=\liminf_{x\to\infty}\frac{P(S_{(\tau)}>x)}{\overline{F}(x)}
=E\tau. \label{randomly-eq5}
\end{equation}
Moreover, if $F\in\cal{S}_*$, then
\begin{equation}
\lim_{x\to\infty}\frac{P(S_{\tau}>x)}{\overline{F}(x)}=\lim_{x\to\infty}
\frac{P(S_{(\tau)}>x)}{\overline{F}(x)}=E\tau.\label{randomly-eq6}
\end{equation}

\item[\rm {(ii)}] Assume that $EX_1\ge 0$ and that there exists $c>E\xi_1$ such that $P(c\tau>x)=o(\overline{F}(x))$ as $x\to\infty$. If $F\in\cal{S}$, then
\begin{equation}
\liminf_{x\to\infty}\frac{P(S_{\tau}>x)}{\overline{F}(x)}=\liminf_{x\to\infty}
\frac{P(S_{(\tau)}>x)}{\overline{F}(x)}=E\tau. \label{randomly-eq7}
\end{equation}
Moreover, if $F\in\cal{S}^*$, then
\begin{equation}
\lim_{x\to\infty}\frac{P(S_{\tau}>x)}{\overline{F}(x)}=\lim_{x\to\infty}
\frac{P(S_{(\tau)}>x)}{\overline{F}(x)}=E\tau.\label{randomly-eq8}
\end{equation}

\end{enumerate}
\end{theorem}

\medskip
\noindent {\bf Proof.} \ We first prove (i). Since $F\in\cal{S}$ is
heavy tailed and has a finite mean, it follows from Lemma~4 of Foss
and Korshunov [11] that there exists a sequence $\{x_{l}\}$ such
that $x_{l}\to\infty$ as $l\to\infty$ and
\begin{equation}
\lim_{l\to\infty}\frac{1}{\overline{F}(x_{l})}\int_0^{x_{l}}
\overline{F}(x_{l}-y)\overline{F}(y)dy=2\int_0^{\infty}\overline{F}(z)dz.
\label{randomly-eq9}
\end{equation}
Lemma 9 of Denisov et al. [5] implies that
$\overline{F_h*F_h}(x_{l})\sim 2\overline{F_h}(x_{l})$ uniformly in
$h\in [1,\infty)$ as $l\to\infty$, where $F_h$ is defined in
(\ref{ee}). If $EX_1<0$, it follows from the result of Korshunov
[19]   that
\begin{equation}
P(S_{(n)}^*>x_{l})\sim
\frac{1}{|E{X_1}|}\int_{x_{l}}^{x_{l}+n|E{X_1}|}\overline{F}(y)dy,\label{randomly-eq10}
\end{equation}
uniformly in $n\ge 1$. Consider the relation
\begin{equation}
P(S_{(\tau)}>x_l)\sim E\tau \overline{F}(x_l).\label{randomly-eq11}
\end{equation}
Since $c(u_1,\cdots,u_n)<M$ for all $(u_1,\cdots,u_n)\in [0,1]^n$,
then for any $x$,
\begin{equation}
P(S_{(n)}>x)\le M P(S_{(n)}^*>x).\label{randomly-eq12}
\end{equation}
Thus, from (\ref{randomly-eq10}), we have
\begin{equation}
P(S_{(n)}>x_{l})\le
M(1+o(1))n|E{X_1}|\overline{F}(x_l), \nonumber %%\label{randomly-eq13}
\end{equation}
for all $n\ge 1$.  Applying the dominated convergence theorem and
(\ref{finitely-eq4}), we obtain
\begin{equation}
\lim_{l\to\infty}\frac{P(S_{(\tau)}>x_l)}{\overline{F}(x_l)}
=\sum_{k=1}^{\infty}\left(\lim_{l\to\infty}
\frac{P(S_{(n)}>x_l)}{\overline{F}(x_l)}\right)
P(\tau=n)=E\tau,\label{randomly-eq14}
\end{equation}
which proves (\ref{randomly-eq11}). Furthermore, Fatou's lemma gives
\begin{equation}
\liminf_{x\to\infty}\frac{P(S_{\tau}>x)}{\overline{F}(x)}\ge
E\tau,\label{randomly-eq15}
\end{equation} without any restriction on the sign of $EX_1$.  Since
$P(S_{\tau}>x)\le P(S_{(\tau)}>x)$ for all $x$, (\ref{randomly-eq5})
follows from (\ref{randomly-eq14}) and (\ref{randomly-eq15}). On the
other hand, if $F\in\cal{S}_*$ and $EX_1<0$, then
(\ref{randomly-eq10}) holds for any $x$, so that
(\ref{randomly-eq6}) can be proved by modifying the proof of
(\ref{randomly-eq5}).

To prove (ii), it is sufficient to prove (\ref{randomly-eq11})
%%$$P(S_{(\tau)}>x_l)\sim E\tau \overline{F}(x_l),$$
for some sequence $\{x_l\}$. Since $F\in\cal{S}$, it follows from
(\ref{finitely-eq4}) that$$P(S_{(n)}>x)\sim n\overline{F}(x).$$
Thus, there exists an increasing function $N(x)\to\infty$ such that
$$P_1(x):=P(S_{(\tau)}>x, \tau\le N(x))\sim E\tau \overline{F}(x).$$
Let $\varepsilon=(c-EX_1)/2>0$ and $b=(EX_1+c)/2$. Put
$\tilde{X_i}=X_i-b$ and $\tilde
{S_n}=\tilde{X_1}+\cdots+\tilde{X_n}$. Then,
$E\tilde{X_i}=-\varepsilon<0$. Using (\ref{randomly-eq10}) and
(\ref{randomly-eq12}), we have
\[ P(S_{(n)}>x_l)\le M
P(S_{(n)}^*>x_l)\le M P(\tilde{S}_{(n)}^*>x_l-bn).
\]
Following the steps of the proof of Theorem 1 (ii) of Denisov et al.
[8], one gets
$$P_2(x_l):=P(S_{(\tau)}>x_l, \tau\in(N(x_l),x_l/c]) =o(\overline{F}(x_l)).$$
 Finally, the condition $P(c\tau>x)=o(\overline{F}(x))$ gives
$$P_3(x):=P(S_{(\tau)}>x, c\tau>x)=o(\overline{F}(x)).$$
Thus, $$P(S_{(\tau)}>x_l)\equiv P_1(x_l)+P_2(x_l)+P_3(x_l)\sim E\tau
\overline{F}(x_l)\ \ as \ \ l\to\infty.$$ Also, if $F\in\cal{S}^*$,
then $$P(S_{(\tau)}>x)\equiv P_1(x)+P_2(x)+P_3(x)\sim E\tau
\overline{F}(x), \ \ as \ \ x\to\infty.$$ Hence, the proof of (ii)
is complete. \hfill $\Box$

\medskip
\noindent {\bf Remark 4.1.} \ Theorem 4.4 partially extends
(\ref{finitely-eq3}) and (\ref{finitely-eq4}) to the case of random
sums. Under the independence assumption and for $F\in \cal{S}^*$,
(\ref{randomly-eq6}) was established in Denisov et al. ([8, Theorem
1(i)]).  Here, (\ref{randomly-eq8}) generalizes the result of
Denisov et al. ([8, Theorem 1(ii)]) to the dependent case. For
related work, we refer the reader to Ng et al. [22], and Ng and Tang
[21]. Concerning the asymptotics for the maximum of $S_{(\tau)}$, it
was shown in Foss and Zachary [12] that the relation
$P(S_{(\tau)}>x)\sim E\tau \overline{F}(x)$ holds for any stopping
time $\tau\le \infty$ and $F\in\mathcal{S}_*$. \hfill $\Box$

%\vskip 0.2cm
%\noindent{\bf 5. Applications to risk theory}
%\vskip 0.2cm

\setcounter{equation}{0}
\section{Applications to risk theory}\label{risk}

In this section, we present an example to illustrate some
applications of our main results.

\medskip
\noindent {\bf Example 5.1.} \ Following the formulation of Ng et
al. [20], we can write the surplus of an insurance company at time
$t$ as
\begin{equation}
U_{\delta}(t)=x e^{\delta t}-\sum_{k=1}^{N(t)}X_k e^{\delta
(t-\sigma_k)},\; t\ge 0, \label{risk-eq1}
\end{equation}
where $x\ge 0$ is the initial surplus, $\delta\ge 0$ is the constant
interest force, $\sigma_k$ is the time at which the $k$th customer
arrives and buys an insurance contract with $\sigma_0=0$,
$N(t)=\max\{k\ge 0:\sigma_k\le t\}$ is the individual
customer-arrival process, and $X_k$ represents the total potential
claims due to the $k$th customer. In this example, $\{X_k, k \ge
1\}$ is a sequence of random variables which are not necessarily
i.i.d.  Furthermore, if $\sigma_k=k$ for each $k\ge 1$, then risk
process (\ref{risk-eq1}) can be rewritten as
\begin{equation}
 U_{\delta}(0)=x,\;\;  U_{\delta}(n)=x(1+r)^n-\sum_{k=1}^n X_k (1+r)^{n-k}, \quad n=1,2.\cdots, \label{risk-eq2}
\end{equation}
where $r=e^{\delta}-1$. In the literature, model (\ref{risk-eq2})
corresponds to a discrete-time insurance risk model with a constant
interest rate (see, for example, Tang [24]).

Let $V_k= X_k (1+r)^{-k}$ and $S_n=\sum_{k=1}^n V_k.$ Then, we can
rewrite (\ref{risk-eq2}) as $U_{\delta}(n)=(1+r)^n(x- S_n)$. Define
$n$-period finite-time ruin probability as
\begin{equation}
\psi_n(x)=P\left(\inf_{1\le k\le n}  U_{\delta}(k) < 0|
 U_{\delta}(0)=x\right)=P\left(\max _{1\le k\le n} S_k>x\right). \nonumber %%\label{risk-eq3}
\end{equation}
From Theorem~\ref{thrm3-1}, we obtain

\begin{corollary}\label{corollary5-1} \
Assume that $X_k \sim F_k$ and hence $V_k$ are heavy-tailed random
variables for $k= 1,2,\cdots$. Under assumption {\bf H3}, if  $F_1,
F_2,\cdots, F_n\in \cal{L}$, then
\begin{equation}
\liminf_{x\to\infty}\frac{\psi_n(x)}{\sum_{k=1}^n
\overline{F_k}(x(1+r)^k)}=1; \nonumber %%\label{risk-eq4}
\end{equation}
and if additionally  $P(\sum_{k=1}^n
V_k^{*+}>x)\sim\sum_{k=1}^{n}\overline{F_k}(x)$,
 then
\begin{equation}
\lim_{x\to\infty}\frac{\psi_n(x)}{\sum_{k=1}^n
\overline{F_k}(x(1+r)^k)}=1.\label{risk-eq5}
\end{equation}
\end{corollary} \hfill $\Box$

Note that if $\delta=0$ and $ X_k=Z_k-(1+\rho)\mu,$  where  $\{Z_k,
k \ge 1\}$ is a sequence of nonnegative random variables with common
distribution function F and finite mean $\mu$, and the positive
constant $\rho$ can be interpreted as the safety loading, then
(\ref{risk-eq1}) has the form
\begin{equation}
U_{\delta}(t)=x-\sum_{k=1}^{N(t)}X_k\equiv
x-S_{N(t)}. \nonumber %%\label{risk-eq6}
\end{equation}
which is the so-called customer-arrival-based insurance risk model
studied by Ng et al. [20].

For the ruin probability within a finite horizon $T$ given by
\[
\psi (x; T)=P\left(\inf_{0\le t\le T}U(t)<0\right)=P(S_{(N(T))}>x),
\]
one can apply Theorem~\ref{thrm4-4} to obtain
\begin{corollary}\label{corollary5-2} \
Let $\{Z_k, k\ge 1\}$ independent of $N(t)$ be a sequence of
heavy-tailed random variables with common distribution $F$.  Under
the conditions in Theorem~\ref{thrm4-4}, if $F\in\cal{S}$, then
\[
\liminf_{x\to\infty}\frac{\psi (x; T)}{\overline{F}(x)}=\lambda(T).
\]
Also, if $F\in\cal{S}^*$, then
\begin{equation}
\lim_{x\to\infty}\frac{\psi (x;
T)}{\overline{F}(x)}=\lambda(T). \nonumber %%\label{risk-eq7}
\end{equation}
\end{corollary} \hfill $\Box$

\medskip
\noindent {\bf Remark 5.1.} \   Under the independent setting of
$\{Z_k, k\ge 1\}$ and other conditions, (\ref{risk-eq5}) were
obtained in several papers. For example, Ng et al. [22] obtained
(\ref{risk-eq5}) under the conditions that
$F\in\mathcal{L}\cap\mathcal{D}$ and $P(N(T)>x)=o(\overline{F}(x))$;
Ng et al. [20] obtained (\ref{risk-eq5}) under the conditions that
$\{Z_k, k\ge 1\}$ is a sequence of i.i.d. subexponential random
variables and $E(1+\varepsilon)^{N(T)}<\infty$ for some
$\varepsilon>0$; and Kass and Tang [14] weakened the condition on
$N(\cdot)$ and obtained (\ref{risk-eq5}) under the condition that
$F$ is strongly subexponential, that is, $F\in \mathcal{S}_*$.
\hfill $\Box$

\noindent{\bf Acknowledgements.} \ The authors are thankful to
Professor Qihe Tang
 for his helpful comments and constructive
suggestions, which have considerably enhanced this work.   The
research of Kam C. Yuen was supported by a university research grant
of the University of Hong Kong. The research of Chuancun Yin was
supported by the National Natural Science Foundation of China (No.
10771119) and the Research Fund for the Doctoral Program of Higher
Education of China (No. 20093705110002).

\end{document}